\documentclass[10pt]{article}
\usepackage{amssymb}
\usepackage{amsfonts}
\usepackage{amsmath}
\usepackage{graphicx}
\usepackage{amsfonts}
\usepackage{amssymb}
\usepackage{epstopdf}
\usepackage{color}

\usepackage[section]{placeins}

\newcommand{\la}{\lambda}

\newcommand{\cf}{{\cal F}}

\newcommand{\cac}{\mathcal C}

\newcommand{\ka}{\kappa}

\newcommand{\al}{\alpha}
\newcommand{\ga}{\gamma}
\newcommand{\ep}{\varepsilon}
\newcommand{\si}{\sigma}

\newcommand{\R}{{\mathbb R}}

\newcommand{\lp}{\left(}
\newcommand{\rp}{\right)}
\newcommand{\lc}{\left[}
\newcommand{\rc}{\right]}

\newtheorem{theorem}{Theorem}[section]

\newtheorem{hypothesis}[theorem]{Hypothesis}

\newtheorem{proposition}[theorem]{Proposition}

\begin{document}
\thispagestyle{empty}

\begin{center}

\huge { Coinfection in a stochastic model for bacteriophage systems}

\vspace{.5cm}

\normalsize {\bf Xavier Bardina\footnote{X. Bardina is partially
supported by the grant MTM2015-67802-P from MINECO.}, Carles
Rovira$^{*}$\footnote{C. Rovira is partially supported by the grant
MTM2015-65092-P from MINECO/FEDER, UE}}

{\footnotesize \it $^1$ Departament de Matem\`atiques,
Universitat Aut\`onoma de Barcelona,
08193 Bellaterra.

$^3$ Departament de Matem\`atiques i Inform\`atica, Universitat de Barcelona, Gran
Via 585, 08007-Barcelona.

 {\it E-mail addresses}: Xavier.Bardina@uab.cat,
Carles.Rovira@ub.edu}

{$^{*}$corresponding author}

\end{center}

\begin{abstract}
A system modeling bacteriophage treatments with
coinfections in a noisy context is analyzed. We prove that in a small noise
regime, the system converges in the long term to a bacteria free
equilibrium. Moreover, we compare the treatment with coinfection
with the treatment without coinfection, showing how the coinfection
affects the  dose of bacteriophages that is needed to eliminate the
bacteria and the velocity of convergence to the free bacteria
equilibrium.
\end{abstract}



{\bf AMS 2000 MSC:}  92D25

{\bf Running head:} Coinfection in a stochastic model

\section{Introduction}

The emergence of pathogenic bacteria resistant to most currently available antimicrobial agents has become a critical problem in  medicine.
 The development of alternative antiinfection modalities has become a priority. Bacteriophage therapies are one of these  alternatives.
Prior to the discovery and widespread use of antibiotics, it has been
suggested that bacterial infections could be treated by the
administration of bacteriophages, but early clinical studies
with bacteriophages were not pursued in the United States and
Western Europe. Nowadays, these therapies are reemerging and
attracting the attention of the scientific community.

Let us explain the (lytic) bacteriophage mechanism: after
attachment, the virus' genetic material penetrates into the bacteria
and uses the host's replication  mechanism to self-replicate. After
some time $\tau$, the bacteria encounters death releasing some new
viruses, ready to attack other bacteria.

When a bacteria has absorbed a virus particle it  is called an infection of a bacteria by a bacteriophage (virus).  An infected bacteria,
 that is, one that has already been previously absorbed and infected, can have a second or additional adsorption. Some authors have used the term superinfection
 as a synonym for what we have described  as secondary adsorption. On the other hand, the word coinfection denotes, generally, the infection of a single cell
  by more than one virus. Particularly, it means the infection of a single bacteria by two or more bacteriophages of the same type. Coinfection is similar but not
  identical to superinfection, since coinfection does not necessarily follows multiple adsorptions of a single bacteria since later phages can be blocked from
  entering the bacteria.

There is a long history of mathematical modelling of  phage dynamics. One
of the first papers was the work of  Campbell \cite{C} where he
proposed a model based on a system  of delay differential equations.
Deterministic models can be found, for instance, in    \cite{BL},
\cite{CTJ}, \cite{CPR}, \cite{LB}, \cite{LSC}, \cite{PJ},
\cite{WBH}. The literature about stochastic models is scarce. For
instance, in \cite{ST} the authors give a stochastic model allowing
multiple bacteriophage adsorption to host. On the other hand,
\cite{BS} was one of the first papers dealing with co/infection and
superfinfection models in evolutionary epidemiology, as previous
models  took only first infections into account.  A general
discussion about how superinfections and coinfections have been
modeled in evolutionary epidemiology can be found in \cite{Al}.

\bigskip

In \cite{BBRT} we have considered   a stochastic model  with a  constant
injection of phages into the system. This variant corresponds to a
treatment for cattle against Salmonella, which was  brought to
our attention by the Molecular Biology Group of the Department of
Genetics and Microbiology at \emph{Universitat Aut\`onoma de
Barcelona}. We  modeled a  bacteriophage system to a kind of
predator-prey equation.

Set $S(t)$ (resp. $Q(t)$) for the non-infected bacteria (resp.
bacteriophages) concentration at time $t$. Consider  a truncated
identity function $\sigma:\R_+\to\R_+$, such that
$\si\in\cac^{\infty}$, $\si(x)=x$ whenever $0\le x\le M$ and
$\si(x)=M+1$ for $x>M+1$. Then the model  is as follows:
\begin{equation}\label{e:truncada2}
\left\{
\begin{aligned}
dS(t) & =\left[\alpha-k_1\sigma(Q(t))\right]S(t)dt\\
dQ(t )& = \left[d-mQ(t)-k_1 \sigma(Q(t))S(t)+k_1 \, b \, e^{-\mu \tau} \sigma(Q(t-\tau))S(t-\tau)\right] dt,
\end{aligned}
\right.
\end{equation}
where $\alpha$ is the reproducing rate of the bacteria, $k_1$ is the
adsorption rate, $d$ also stands for the quantity  of bacteriophages
inoculated per unit of time, $m$ is their death rate,  $b$ is the
number of bacteriophages which are released after replication within
the bacteria cell, $\tau$ is the delay necessary for the reproduction
of bacteriophages (called latency time) and the coefficient $e^{-\mu
\tau}$ represents an attenuation in the release of bacteriophages
(given by the expected number of bacteria cell's deaths during the
latency time, where $\mu$ is the bacteria's death rate).
In fact, $\alpha=\beta - \mu$ where $\beta$ is the bacteria's reproduction rate.
A given
initial condition $\{S_{0}(t), Q_{0}(t); -\tau\le t\le 0\}$ is also
specified.

Given a large enough $M$ we show that when $k_1d/m>\al$, there exists
a unique stable steady state $E_0=(0,d/m)$  (bacteria have been
eradicated), while when $k_1d/m< \al$, the point $E_0$ is still an
equilibrium but it becomes unstable and there exists another coexistence
equilibrium. This paper only studies results regarding the bacteria-free
equilibrium $E_0$, since it corresponds to the main practical
situation, where high doses  of phages are usually introduced in 
cattle feed. Actually we also consider a small random perturbation
of the form
\begin{equation}\label{e:truncada-delayed}
\left\{
\begin{aligned}
d S^\varepsilon (t) & =\left[\alpha-k_1\sigma(Q^\varepsilon(t))\right]S^\varepsilon(t)dt + \varepsilon\sigma(S^\varepsilon(t))\circ dW^1(t)\\
d Q^\varepsilon(t) & = \left[d-m Q^\varepsilon(t)-k_1 \sigma(Q^\varepsilon(t))S^\varepsilon(t)+k_1 \, b \, e^{-\mu \tau} \sigma(Q^\varepsilon(t-\tau))S^\varepsilon(t-\tau)\right] dt \\
& \qquad + \varepsilon\sigma(Q^\varepsilon(t))\circ d W^2(t),
\end{aligned}
\right.
\end{equation}
where $\varepsilon$ is a small positive coefficient and $W=(W^1,W^2)$ is a
2-dimensional Brownian motion and with Stratonovich differentials, denoted by $\circ\, dW$. We get a concentration result for the perturbed system around $E_0$.

\bigskip

Our aim is to study the problem of coinfection in the model we have
presented in \cite{BBRT}. Due to the ambiguity in terminology we
have to specify what we understand by coinfection. After the first
infection by a bacteriophage and before the death of the bacteria
(we assumed that there  is a period of  time $\tau$) when the bacteriophages can infect
 the bacteria more times. These later adsorptions will not  affect
the behaviour of the bacteria but they cause the destruction of the
attacker bacteriophages. Thus, we will loose some bacteriophages.

In order to describe this situation, we introduce a new process $I(t)$, that gives the infected bacteria concentration at time $t$.
Thus, we transform the model  (\ref{e:truncada2}) into the next one
\begin{equation}
\label{MOD2}
\left\{
\begin{array}{lll}
dS(t) & = & (\alpha - k_1 \sigma(Q(t))) S(t) dt
\\
dI(t) & = & [k_1 \sigma(Q(t)) S(t)  - \mu I(t)  - k_1 e^{-\mu \tau }
\sigma(Q(t-\tau)) S(t-\tau)]dt
\\
dQ(t) & = & [d-m Q(t) -k_1 \sigma(Q(t)) S(t)- k_2 \sigma(Q(t)) I(t)
\\ && + k_1 b  e^{-\mu \tau } \sigma(Q(t-\tau)) S(t-\tau)]dt.
\end{array}
\right.
\end{equation}
Notice that the variation on $I(t)$ increases with the number of
infected bacteria $ k_1 \sigma(Q(t)) S(t)$ and decreases with the rate of
dead  infected bacteria $\mu I(t)$ and the infected bacteria that
disappear after the latency time. On the other hand, in $dQ(t)$
 the term $-k_2 \sigma(Q(t)) I(t)$ appears that explains the bacteriophages
that the system looses when they try to infect infected bacteria
(coinfection). Now $k_1$ is the adsorption rate for noninfected bacteria and $k_2$ is the adsorption rate by infected ones. Moreover, when $k_2=0$ we get the system without
coinfection. Thus, our  model (\ref{MOD2}) is an extension of the model
presented in \cite{BBRT}.

We will show that under certain conditions
  a unique stable steady state $E_0=(0,0,d/m)$  (free-bacteria equilibrium) exists  and we obtain
  also a concentration result around $E_0$ for a perturbed system. These results are analogous to those obtained in \cite{BBRT}.  Furthermore, we will
  compare both models to determine the role of the coinfection in the behaviour of the system.

\bigskip

Our article is structured as follows: Section \ref{main} is devoted
to giving the main results with the corresponding biological discussion.
Then we show the detailed mathematical analysis in Section
\ref{maths}.

\section{Main results and discussion}\label{main}
In this section we will study our coinfection model (\ref{MOD2}) and the fluctuations of its corresponding stochastic model. Before going on with the study of the deterministic model, let us present the set of hypothesis on the coefficient $\sigma$ and on the initial condition. The hypothesis on $\sigma$ will be the same than those in \cite{BBRT}.

\begin{hypothesis}\label{hyp:sigma}
The coefficients of our differential systems satisfy the following assumptions:

\smallskip

\noindent
\emph{(i)}
The function $\sigma:\R_+\to\R_+$ is such that $\sigma\in\cac^{\infty}$, and satisfies $\sigma(x)=x$ for $0\le x\le M$ and $\sigma(x)=M+1$ for $x>M+1$. We also assume that $0\le \sigma'(x)\le C$ for all $x\in\R_+$, with a constant $C$  such that $C>1$.

\smallskip

\noindent
\emph{(ii)}
As far as the initial condition is concerned, we assume that it is given as continuous positive functions $\{S_{0}(s), Q_{0}(s); -\tau\le s\le 0\}$ and a constant $I_0$.
\end{hypothesis}

\subsection{Deterministic model}

We need to introduce a new assumption on the initial condition in order to get the global existence and positivity of solution.

\begin{hypothesis}\label{condi}
Assume that the initial condition satisfies:
$$I_0 \ge\ k_1 e^{-\mu \tau } \int_{-\tau}^ 0 \sigma(Q(s)) S(s) ds.$$
\end{hypothesis}

Clearly, if $Q(t)=0$ for any $t \in [-\tau,0]$ then the condition is only $I_0\ge0$. It corresponds to the case when we begin to introduce bacteriophages at time $t=0$.

Under Hypothesis  \ref{hyp:sigma} and \ref{condi} we have  the positivity and the global existence of solution (see Proposition \ref{exis_pos}). Nevertheless, to obtain the existence and stability of equilibrium points we need more hypothesis: a set of hypothesis on the initial condition and a set of hypothesis on the coefficients.

\begin{hypothesis}\label{hyp:coeff-delay}
    We will suppose that the coefficients satisfy the following conditions, valid for  any $t\in[-\tau,0]$:

    \smallskip

    \noindent
    \emph{(i)} The initial condition $(S_{0}(t),I_{0}(t),Q_0(t))$ of the system lies into the region
    $$
    R_0:=\left[0,M\right]\times\left[0,M\right]\times\lc \frac{d(be^{-\mu \tau} \mu)}{mbe^{-\mu \tau} \mu+k_2(mM-d)},M\rc.
    $$

    \smallskip

    \noindent
    \emph{(ii)}
    We have $(m b\, e^{-\mu\tau} \mu + k_2(mM-d)) Q_{0}(t) S_{0}(t)> d\mu S_{0}(0)$, and $b\, e^{-\mu\tau} >1$.

    \smallskip

    \noindent
    \emph{(iii)}
    The condition $S_{0}(t)<\frac{mM-d}{k_1 be^{-\mu\tau}M}$.

    \noindent
    \emph{(iv)}
    $I_{0}(0)<\frac{mM-d}{be^{-\mu\tau}\mu}$.

\end{hypothesis}

\begin{hypothesis}\label{hyp:coeff-sde}
    We will suppose that the coefficients satisfy that $\frac dm < M$ and  $$d > \frac{\alpha m }{k_1} \frac{be^{-\mu \tau} \mu+k_2(M-\frac dm)}{be^{-\mu \tau} \mu}.$$
\end{hypothesis}
Then under Hypothesis  \ref{hyp:coeff-sde}
$$\frac{\alpha }{k_1} < \frac{d(be^{-\mu \tau} \mu)}{mbe^{-\mu \tau} \mu+k_2(mM-d)}<\frac dm < M.$$

Notice again that when $k_2=0$ we get the same hypothesis as the
model without coinfection. Moreover, when $k_2$ is
increasing, we find that the constant dose $d$ must increase, i.e.,
if we loose more bacteriophages by coinfection we need to introduce
a bigger dose of them. On the other hand, the region where  the
initial condition $Q_0$ lives can have smaller lower boundary. It means that since
the dose will be bigger, the concentration of viruses in the initial
condition can be smaller.

\smallskip

Going back to the study of our model, under Hypothesis  \ref{hyp:sigma}, \ref{condi}, \ref{hyp:coeff-sde} and \ref{hyp:coeff-delay}, we can prove the
 boundedness of the solution of the system (see Proposition \ref{prop:inv-region-delay}) and we are able to study the existence of equilibrium points.

 Let us study the equilibrium points. We have to find the solution to the following equations:
\begin{equation}
\label{MOD2EQUI}
\left\{
\begin{array}{lll}
0 = & (\alpha - k_1 \sigma(Q)) S
\\
0= & k_1 \sigma(Q) S - \mu I -  k_1  e^{-\mu \tau } \sigma(Q) S
\\
0= & d-m Q -k_1 \sigma(Q) S - k_2 \sigma(Q)   I +  k_1 b  e^{-\mu \tau } \sigma(Q) S.
\end{array}
\right.
\end{equation}
Clearly, when $S=0$ we get that $I=0$ and $Q=\frac{d}{m}$. So, the bacteria-free equilibrium $E_0=(0,0,\frac{d}{m})$ exists.
To ensure the existence of this equilibrium point we need that $\frac{d}{m}<M$, that it is true under Hypothesis \ref{hyp:coeff-sde}.
In the case $M+1 < \frac{\alpha}{k_1}$ it is clear that no other equilibrium exists.
Furthermore, if $M \ge \frac{\alpha}{k_1}$ a possible coexistence equilibrium should be
$$Q_c=\frac{\alpha}{k_1}, \qquad I_c(S)=\frac{\alpha}{\mu}(1- e^{- \mu \tau})S,
$$
and $$
S_c=\frac{\mu (k_1 d - m \alpha)}{ \alpha (\mu k_1 (1-  b  e^{-\mu \tau }) +   \alpha  k_2
    (1- e^{- \mu \tau}) )}.$$
More precisely, it will only exist if $d \le m \frac{\alpha}{k_1} $
and  $b e^{-\mu \tau }   > (\frac{\alpha}{\mu} \frac{k_2}{k_1}(1-
e^{-\mu \tau } )+1)$ or if $d \ge  m \frac{\alpha}{k_1} $ and  $1< b
e^{-\mu \tau }    <  (\frac{\alpha}{\mu} \frac{k_2}{k_1} (1-
e^{-\mu \tau } )+1)$. From the biological point of view, these
situations correspond to the cases of a ``small" dose of
``efficient" viruses and  or ``large" dose of ``nonefficient"
viruses, respectively.

\bigskip

As we have explained in the introduction we are interested in the
behaviour of the free disease equilibrium $E_0=(0,0,\frac{d}{m})$.
We get that if   $M > \frac{\alpha}{k_1}$  and $\alpha - k_1
\frac{d}{m} \le 0$,   the bacteria-free equilibrium $E_0$ is the
unique steady state and if $M> \frac{d}{m}$ it is asymptotically
stable (see Proposition \ref{stable}).  We can also get the
exponential convergence to the bacteria-free equilibrium point.

\begin{theorem}\label{thm:exp-cvgce-deterministic}
    Assume Hypothesis \ref{hyp:coeff-sde}, \ref{hyp:sigma}, and \ref{hyp:coeff-delay} are satisfied, and let $R$ be the region defined at Proposition \ref{prop:inv-region-delay}. Then the solution of system \ref{MOD2})
 with initial condition $(S_0,I_0, Q_0)\in R$ exponentially converges to the equilibrium $E_0$:
    \begin{equation}\label{eq:exp-cvgce-E0}
    |(S(t),I(t),Q(t))-E_0| \le c \, e^{-\eta t},
    \quad\mbox{with}\quad
    \eta= \ga \wedge m \wedge \mu,
    \end{equation}
    where  $\ga=\tfrac{k_1d}{m}-\alpha>0$.
\end{theorem}

Summarizing, the free equilibrium point is, in some sense, the same point that in the model without coinfection. That is, the concentration of bacteria is 0 and the concentration of bacteriophages is $\frac dm$. We also have exponential convergence but in our model with coinfection it will be slower  or equal that in the model without coinfection. More precisely, if in  \cite{BBRT} it was of order  $e^{-(\ga \wedge m)  t}$ in the model with coinfection it will be of order $e^{-(\ga \wedge m \wedge \mu) t}$.

\subsection{Stochastic fluctuations}

Let us analyze the stochastic case. We will introduce the random
effects following the ideas we have used in \cite{BBRT}. We will
assume that the noise enters in a bilineal way that ensures
positivity of the solution.  Thus, we consider system (\ref{MOD2})
with a small random perturbation of the form
\begin{equation}
\label{MOD2s}
\left\{
\begin{array}{ll}
dS^\varepsilon(t) = & \big(\alpha - k_1 \sigma(Q^\varepsilon(t))\big) S^\varepsilon(t) dt + \varepsilon \sigma(S^\varepsilon(t)) \circ dW^1(t),
\\
dI^\varepsilon(t) = & \big[ k_1 \sigma(Q^\varepsilon(t)) S^\varepsilon(t)  - \mu I^\varepsilon(t)  - k_1 e^{-\mu \tau } \sigma(Q^\varepsilon(t-\tau)) S^\varepsilon(t-\tau) \big] dt \\
dQ^\varepsilon(t) = & \big[ d-m Q^\varepsilon(t) -k_1 \sigma(Q^\varepsilon(t)) S^\varepsilon(t) - k_2 \sigma(Q^\varepsilon(t))I^\varepsilon(t) \\ & + k_1 b  e^{-\mu \tau } \sigma(Q^\varepsilon(t-\tau)) S(t-\tau) \big]dt
+ \varepsilon \sigma(Q^\varepsilon(t)) \circ dW^2(t)
\end{array}
\right.
\end{equation}
where $\varepsilon$ is a small positive coefficient and $W=(W^1,W^2)$ is a
2-dimensional Brownian motion defined on a complete probability
space $(\Omega,\cf,P)$ equipped with the natural filtration $(\cf_t)_{t\ge 0}$ associated to the Wiener process $W$.  Recall that $\circ dW(t)$ denotes a Stratonovich integral.

The existence follows from the fact that the coefficients of the equation are locally Lipschitz with linear growth (see Theorem 2.7 in \cite{BBRT}).
The positivity holds using the same arguments that in Proposition 2.8 in \cite{BBRT}.

Let us introduce some notation. For a continuous function $f$, we set $\|f\|_{\infty,L}=\sup_{x\in L}|f(x)|$. Set
$Z^\varepsilon=(S^\varepsilon, I^\varepsilon, Q^\varepsilon).$
Then we can state the result about convergence to $E_0$ as follows:

\begin{theorem}\label{thm:concentration-equilibrium}
    Given positive initial conditions and  Hypothesis  \ref{hyp:coeff-sde}, \ref{hyp:sigma}, and \ref{hyp:coeff-delay} , equation (\ref{MOD2s}) admits a unique solution which is almost surely an element of $\cac(\R_+,\R_+^3)$. Set
      $\eta=m\wedge \ga \wedge \mu$ and consider three
     constants $1<\ka_1<\ka_2<\ka_3$. Then there exists $\rho_0$ such that for any $\rho\le \rho_0$ and any interval of time of the form $L=[\ka_1 \ln(c/\rho)/\eta, \ka_2 \ln(c/\rho)/\eta]$, we have
    \begin{equation}\label{eq:concentration-equilibrium}
    P\lp  \|Z^{\varepsilon}-E_0\|_{\infty,L} \ge 2 \rho \rp  \le
    \exp\lp -\frac{c_1 \rho^{2+\la}}{\ep^2}  \rp,
    \end{equation}
    where $\la$ is a constant satisfying $\la>\ka_3/\eta$.
\end{theorem}

    Relation (\ref{eq:concentration-equilibrium}) means that the kind of deviation we might expect from the noisy system (\ref{MOD2s}) with
    respect to the equilibrium $E_0$ is of order $\varepsilon^{2\vartheta }$ with $\vartheta =2\eta/\ka_3$. This range of deviation happens at
    a time scale of order $\ln(\rho^{-1})/\eta$. As in  the exponential convergence for deterministic model, the convergence of the stochastic model
     with coinfection will be slower  or equal that the convergence of the stochastic model without coinfection.

\bigskip

\section{Mathematical analysis and proofs}\label{maths}

In this section we will state the results described in Section \ref{main} and we will prove the Propositions and Theorems presented in the same section.
Since some of the proofs are similar to those given in \cite{BBRT}, we only will give some details of the proofs with new arguments and we will refer to those in \cite{BBRT}
in the other cases.

The fist step to analyse the model is to get the existence of a global nonnegative solution.

\begin{proposition}\label{exis_pos}
Under hypothesis \ref{hyp:sigma} and \ref{condi}
the system (\ref{MOD2}) has an unique global nonnegative solution.
\end{proposition}

\noindent {\bf Proof:} Let us study the positivity of the solution. Clearly
$$S(t) =  S(0) \exp \Big(\alpha - k_1 \sigma(Q(t) \Big) \ge 0.$$
On the other hand, if for some $t_0$ it holds that $Q(t_0)=0$ then $Q^{\prime}(t_0) \ge d >0.$ So, $Q(t) \ge 0$ for all $t$. Finally, we can write
\begin{eqnarray*}
I(t) & \ge & I(0) +\int_0^ t  k_1(1-e^{-\mu \tau} ) \sigma(Q(s))
S(s) ds  - \mu \int_0^ t I(s) ds
\\ &  &\qquad - k_1 e^{-\mu \tau } \int_{-\tau}^ 0 \sigma(Q(s)) S(s) ds.
\end{eqnarray*}
Thus, if for some $t_0$ it holds that $I(t_0)=0$ then under Hypothesisi \ref{condi} we have that $I^{\prime}(t_0)=0$. It yields that
 $I(t) \ge 0$ for all $t$.

In order to get the existence of global solution it is enough to
check that  the local solutions are bounded (see for instance
\cite{DGVLW}). Since ${S}^{\prime}(t)\leq \alpha S(t)$, we get that for
all $t>0$, $S(t)\leq S(0)e^{\alpha t}$. On the other hand,
${Q}^{\prime}(t)\leq d+k_1be^{-\mu\tau}\sigma(Q(t-\tau))S(t-\tau)$. Using
that $\sigma(x)\leq x$ we get that
$${Q}^{\prime}(t)\leq d+ k_1be^{-\mu\tau}S(0)e^{\alpha t}Q(t-\tau).$$
Applying a Gronwall's lemma (see \cite{K} Lemma A.1) we obtain that
$$Q(t)\leq (Q(0)+dt+k_1bS(0) e^{-\mu \tau}\int_{-\tau}^ 0  e^{\alpha s} ds )\exp\left(k_1b S(0) e^{-\mu \tau}  \int_{0}^t e^{\alpha
s}ds\right).$$ Finally, notice  that ${I}^{\prime}(t)\leq
k_1\sigma(Q(t))S(t)\leq k_1Q(t)S(t)$. So, fixed $T$, the local
solutions are bounded in $[0,T]$. \hfill$\square$

\bigskip

For simplicity, set
$$\nu=\frac{d(be^{-\mu \tau} \mu)}{mbe^{-\mu \tau} \mu+k_2(mM-d)}.$$
Now, we can also prove the following proposition that gives us the boundedness of the solution:

\begin{proposition}\label{prop:inv-region-delay}
    Under Hypothesis  \ref{hyp:sigma}, \ref{condi}, \ref{hyp:coeff-sde} and \ref{hyp:coeff-delay}, the region
    \begin{eqnarray*}
    \!\!\!&\!\!\!\!\!\!&\!\!\!
    R:=R_1 \times R_2 \times R_3 \\\!\!\!&\!\!\!=\!\!\!&\!\!\!\left[0,\frac{mM-d}{k_1be^{-\mu\tau}M}\right]\times
    \left[0,\frac{mM-d}{be^{-\mu\tau}\mu}\right]\times\left[\frac{d(be^{-\mu \tau} \mu)}{mbe^{-\mu \tau} \mu+k_2(mM-d)},M\right]\subset[0,M]^3
    \end{eqnarray*}
    is left invariant by equation (\ref{MOD2}).
\end{proposition}

\noindent{\bf Proof:}
We organize the proof in five steps.

\smallskip

\noindent {\it Step 1: While $Q  \ge \nu$  then  $S \in R_1$ and is
nonincreasing.} Since $S$ is positive it is clear that
\begin{equation*}
S^{\prime}({t})\leq 0  \textrm{ whenever }  Q(t)>\frac\alpha k_1 ,
\quad\mbox{and}\quad S^{\prime}(t)\geq 0  \textrm{ whenever }  Q(t)
< \frac\alpha k_1.
\end{equation*}
On the other hand, our system starts from an initial condition
$$Q_{0}(0)\geq \nu \geq\frac{\alpha}{ k_1}.$$ Thus $S$ is non
increasing and remains in $R_1$ as long as $Q \ge \nu$.

\smallskip

\noindent {\it Step 2: There exists a strictly positive
$\varepsilon$ such that $Q(t)>\nu$ for all $t\in(0,\varepsilon)$. }
Notice that here  a  $\varepsilon_0$  exists such that
$I(t)<\frac{mM-d}{be^{-\mu\tau}\mu}$ for all
$t\in(0,\varepsilon_0)$. So,  we have
\begin{eqnarray*}
Q^{\prime} (0) &\ge & d-m \nu -k_1 \nu S_0(0) -k_2  \nu
\frac{mM-d}{be^{-\mu\tau}\mu }
+k_1be^{-\mu\tau}Q_{0}(-\tau)S_{0}(-\tau)
\\ &=& k_1 \lp be^{-\mu\tau}Q_{0}(-\tau)S_{0}(-\tau)-\nu S_{0}(0) \rp >0,
\end{eqnarray*}
where we have used
Hypothesis (ii) of  \ref{hyp:coeff-delay}.

\smallskip

\noindent {\it Step 3: If $S(t)$ is nonincreasing and $I(t)$ remains
in $R_2$ for any $t \le T$ and  $Q(T) = \nu$ then  $Q^{\prime}
(T)>0$.} Let us consider  what happens when $Q(t_0)=\nu.$   We now
introduce the quantity $t_0=\inf\{t>0:\,\,Q(t)=\nu\}$, and notice
that we have
$$
Q^{\prime}(t_{0})=d-m \nu -k_1 \nu S(t_0) -k_2  \nu
\frac{mM-d}{be^{-\mu\tau}\mu }
+k_1be^{-\mu\tau}\sigma(Q(t_0-\tau))S(t_0-\tau).$$ We can now
distinguish two cases:
\begin{enumerate}
\item If $t_0>\tau$, since $S(t)$ is nonincreasing in $[0,t_0]$,
$S(t_{0}-\zeta)\geq S(t_{0})$ and hence
$$Q^{\prime}(t_{0})\geq k_1S(t_{0})\left(be^{-\mu\tau}\sigma(Q(t_{0}-\tau))-\nu\right)>0,$$
due to the fact that $be^{-\mu\tau}>1$, $M>\nu$ and $Q(t_{0}-\zeta)>\nu$.
\item If $t_0\leq \tau$, since $S(t_{0})\leq S_{0}(0)$ we obtain
$$
Q^{\prime}(t_{0}) \geq  k_1 \left(be^{-\mu\tau}Q_{0}(t_{0}-\tau)S_{0}(t_{0}-\tau)-\nu S_{0}(0)\right)>0,
$$
where we have used again Hypothesis (ii) of  \ref{hyp:coeff-delay}.
\end{enumerate}
This discussion allows thus to conclude that $t_0$ cannot be a finite time.
\smallskip

\noindent {\it Step 4: If $S(t)$ is nonincreasing  for any $t \le T$ and  $Q(T) = M$ then  $Q^{\prime} (T)<0$.}
To this aim notice that, whenever $Q_{0}(0)=M$ we have
 $$
Q^{\prime}(0)\leq
d-mM+k_ 1 be^{-\mu\tau}MS_{0}(-\tau)<0,
$$
where we recall that $S_{0}(-\tau)<\frac{mM-d}{k_1be^{-\mu\tau}M}$ according to Hypothesis \ref{hyp:coeff-delay}. This yields the
existence of $\varepsilon>0$ such that $Q(t)<M$ for all $t\in(0,\varepsilon)$.
We now define $t_1=\inf\left\{t>0:\,\,Q(t)=M\right\}$. It is readily checked that
\begin{eqnarray*}
Q^{\prime}(t_1)&\ge &d-mM+k_1be^{-\mu\tau}\sigma(Q(t_{1}-\tau))S(t_{1}-\tau)\\
&=&d-mM+k_1be^{-\mu\tau}MS(t_{1}-\tau),
\end{eqnarray*}
and we can distinguish again two cases:
\begin{enumerate}
\item If $t_1>\tau$, thanks to the fact that $t\mapsto S(t)$ is non-increasing on $[0,t_1]$, we have
$$
Q^{\prime}(t_1)\leq d-mM+k_1be^{-\mu\tau}MS_{0}(0)<0,
$$
since we have assumed that $S_{0}(0)<\frac{mM-d}{k_1be^{-\mu\tau}
M}$.
\item If $t_1\leq \tau$ then
$$
Q^{\prime}(t_1)\leq d-mM+k_1be^{-\mu\tau}MS_{0}(t_{1}-\tau)<0,
$$
thanks to the fact that $S_{0}(t)<\frac{mM-d}{k_1be^{-\mu\tau}M}$
for all $t\in[-\tau,0]$.
\end{enumerate}
 We have thus shown $Q(t)\le M$ for all $t\ge 0$, which finishes the proof.

\noindent {\it Step 5: While $S$ remains in $R_1$ and $Q$ remains in
$R_3$ then  $I$ lives in $R_2$.}
Notice first that
\begin{eqnarray*}
I^{\prime} (t) & \le & k_1 M   \frac{mM-d}{k_1be^{-\mu\tau} M} - \mu I(t) - k_1 e^{-\mu\tau}\sigma(Q(t-\tau))S(t-\tau) \\
 & = &  \frac{mM-d}{be^{-\mu\tau}} - \mu I(t) - k_1 e^{-\mu\tau}\sigma(Q(t-\tau))S(t-\tau) .
\end{eqnarray*}
We have seen before that $I$ is always nonnegative. Assume now that there exist $t_1$ such that $I(t_1)= \frac{mM-d}{be^{-\mu\tau} \mu}.$
 Then
$$
I^{\prime} (t_1)  \le  - k_1 e^{-\mu\tau}\sigma(Q(t-\tau))S(t-\tau)
\le 0,$$ and so, $I(t)\le \frac{mM-d}{be^{-\mu\tau} \mu}$ for all
$t\ge 0$.

\noindent {\it Conclusion:}
From the previous steps we get that there exists $\varepsilon >0$ such that $(S(t), I(t), Q(t) ) \in R$   for any $t \in [-\tau,\varepsilon)$.
Combining all the steps it is clear that they can not leave the region.
\hfill$\square$

\bigskip

 We can state the stability result as follows:

 \begin{proposition}\label{stable}
    If either $M+1< \frac{\alpha}{k_1}$ or  $M > \frac{\alpha}{k_1}$  and $\alpha - k_1 \frac{d}{m} \le 0$,  system (\ref{MOD2}) has a
    unique steady state $E_0=(0,0,\frac{d}{m})$ . Moreover, the bacteria-free equilibrium $E_0$ is asymptotically stable
     for $\alpha - k_1 \frac{d}{m} <0$ and $M> \frac{d}{m}$.
 \end{proposition}

 \noindent {\bf Proof:}
 When $\tau=0$, using that  $\sigma(Q(t))=Q(t)$, around $E_0$ the differential matrix is:
$$ \left( \begin{array}{ccc}
\alpha - k_1 \frac{d}{m} & 0 & 0\\
0 & -\mu & 0 \\
k_1 (b - 1)\frac{d}{m} &   -k_2 \frac{d}{m}    &-m\end{array} \right),$$
 with  eigenvalues $\lambda_0=\alpha - k_1 \frac{d}{m} $ , $\lambda_1=-\mu<0$ and $\lambda_2=-m<0$. Thus $E_0$ is stable if and only if
$\alpha - k_1 \frac{d}{m} <0$. In order to study the system with
delay, we linearize it around $E_0$, i.e, $S(t)=0+s(t), I(t)=0+i(t)$
and $Q(t)=q(t)+\frac{d}{m}$ and we assume that the solutions are
exponential, i.e.  $s(t)=e^{\lambda t}s, i(t)=e^{\lambda t}i$ and
$q(t)=e^{\lambda t}q.$ We get
\begin{equation}
\label{MOD2_exp}
\left\{
\begin{array}{ll}
\lambda e^{\lambda t}s =  (\alpha - k_1 \frac{d}{m}) e^{\lambda t}s
\\
\lambda e^{\lambda t}i =   k_1 \frac{d}{m}  e^{\lambda t}s   - \mu e^{\lambda t}i  - k_1 e^{-\mu \tau }\frac{d}{m} e^{\lambda ( t-\tau)}s
\\
\lambda e^{\lambda t}q = -me^{\lambda t}q -k_1  \frac{d}{m}  e^{\lambda t}s   -k_2  \frac{d}{m}e^{\lambda t}i + k_1 b  e^{-\mu \tau } \frac{d}{m}e^{\lambda ( t-\tau)}s.
\end{array}
\right.
\end{equation}
Thus  the characteristic equation  is
$$ p(\lambda)= \left| \begin{array}{ccc}
\lambda- (\alpha - k_1 \frac{d}{m})  & 0 & 0\\
 - k_1 \frac{d}{m}    +  k_1 e^{-(\mu+\lambda) \tau }\frac{d}{m} & \lambda+\mu  & 0 \\
k_1  \frac{d}{m} - k_1 b  e^{-(\mu + \lambda) \tau } \frac{d}{m} &   k_2  \frac{d}{m}  &\lambda +m\end{array} \right|=0,$$
and the eigenvalues will be $\lambda_1=\alpha - k_1 \frac{d}{m}<0$, $\lambda_2=-\mu<0$ and $\lambda_3=-m<0$ and so $E_0$ is stable
under the same condition that when $\tau=0$.
\hfill$\square$

\medskip

We can now prove the exponential convergence to the bacteria-free equilibrium point.

\medskip

\noindent{\bf Proof of Theorem \ref{thm:exp-cvgce-deterministic}:} According to Proposition \ref{prop:inv-region-delay}, we have $Q(t)\leq M$ for all $t$.  Doing  now the change of variables $\tilde Q=Q-\frac dm$ we get:
$$\begin{array}{l}
dS(t)=- \lp \gamma S(t) + k_1\tilde Q(t)S(t )\rp \, dt,\\
dI(t)= \lp k_1\frac dm S(t)+k_1\tilde Q(t)S(t)   - \mu I(t) - k_1\frac
dm e^{-\mu\tau}S(t-\tau)  \right. \\
\qquad\qquad  \left.-k_1 e^{-\mu\tau}\tilde Q(t-\tau )S(t-\tau) \rp \, dt,\\
d {\tilde Q}(t)=\lp -m\tilde Q(t)-k_1\frac dm S(t)-k_1\tilde Q(t)S(t)-k_2\frac dm I(t)-k_2\tilde Q(t)I(t)\right. \\
\qquad\qquad  \left. +k_1\frac
dm be^{-\mu\tau}S(t-\tau) +k_1 be^{-\mu\tau}\tilde Q(t-\tau )S(t-\tau) \rp \, dt.
\end{array}
$$
 With our change of variables, we have also shifted our equilibrium to the point $(0,0,0)$. We now wish to prove that $S(t), I(t)$ and $\tilde
Q(t)$ exponentially converge to 0.

\smallskip

The bound on $S(t)$ is easily obtained: just note that
$$
S^{\prime}(t) \leq-\gamma S(t) \, dt,
$$
which yields $S(t)\leq S_{0}(0)\, e^{-\gamma t}$.

As far as $\tilde Q(t)$ is concerned, one gets the bound
\begin{eqnarray*}
\tilde Q^{\prime}(t)&\leq& -m\tilde Q(t)+k_1be^{-\mu\tau}(\frac dm +\tilde Q(t-\tau))S_{0}(0)\, e^{-\gamma(t-\tau)}\\
&\leq& -m\tilde Q(t )+ c\, e^{-\gamma t},
\end{eqnarray*}
with $c=2k_1bMS_{0}(0)\, e^{(\gamma-\mu)\tau}$, and where we have
used the fact that $Q(t)\le M$ uniformly in $t$. Using that equation
$x^{\prime}(t)=-mx(t)+c \, e^{-\gamma t}$ with initial condition
$x_0=\tilde Q_{0}(0)$ can be explicitly solved as
$$
x(t)=\left(\tilde Q_{0}(0) -\frac c{m-\gamma}\right)e^{-mt}+\frac
c{m-\gamma}e^{-\gamma t}
$$
and by comparison, this entails the inequality $\tilde Q(t)\leq c_1\, e^{-(m\wedge \gamma) t}$, where $c_1>0$.

Finally, let us consider $I(t)$. Clearly
$$
I^{\prime}(t) \le  k_1\frac dm S(t)+k_1\tilde Q(t)S(t)   - \mu I(t)
\le  2k_1 M S_0(0) e^{-\gamma t} - \mu I(t).
$$
Following the same method, we get that
$I(t)\leq c_2\, e^{-(\mu \wedge \gamma) t}$.
\hfill$\square$

\bigskip

Finally, we can prove the stochastic convergence.

\bigskip

\noindent {\bf Proof of Theorem \ref{thm:concentration-equilibrium}:} Since we have exponential convergence for the deterministic delayed system, it is enough to check (see subsection 3.1 in \cite{BBRT}) that
 for any $\varepsilon \le \varepsilon(M,T)$
$$
P( \Vert Z^\varepsilon - Z^0 \Vert_{\infty,[0,T]} \ge \rho ) \le
\exp - \big(  \frac{c_2 \rho^2}{e^{K_2 T} \varepsilon^2} \big).
$$
Recall that $\Vert S^0 \Vert_\infty + \Vert I^0 \Vert_\infty\le c_4$ and set $J_1(t)=\int_0^ t \sigma(S^\varepsilon(t)) \circ dW^1(t)$  and
$J_2(t)=\int_0^ t \sigma(Q^\varepsilon(t)) \circ dW^2(t).$ Then

\begin{eqnarray*}
& &\vert S^\varepsilon(t) - S^0(t)  \vert \le \int_0^ t   \vert (\alpha - k_1 \sigma(Q^\varepsilon(s))) (S^\varepsilon(s) - S^0(s)) \vert ds \\ && \qquad
+ \int_0^ t   \vert k_1 (\sigma(Q^\varepsilon(s)) - \sigma(Q^0(s))) S^0(s) \vert ds
+ \varepsilon \vert J_1(t) \vert,
\\
& &\quad\le \int_0^ t   (\alpha + k_1 M)  \vert S^\varepsilon(s) - S^0(s) \vert ds + \int_0^ t   k_1 c_4 C  \vert Q^\varepsilon(s) - Q^0(s) \vert ds
\\ && \qquad
 + \varepsilon \vert J_1(t) \vert,
\\& &\vert Q^\varepsilon(s) - Q^0(s) \vert ds  \le\int_0^ t    m\vert  Q^\varepsilon(s) - Q^0(s) \vert ds  \\
& & \qquad  +  \int_0^ t k_1 b  e^{-\mu \tau } [\vert \sigma(Q^\varepsilon(s-\tau))-\sigma(Q^0(s-\tau)) \vert  \vert S^0(s-\tau) \vert
 \\ & &\qquad \qquad \qquad+\vert S^\varepsilon(s-\tau)-S^0(s-\tau) \vert \vert  \sigma(Q^\varepsilon(s-\tau) \vert ]ds \\
& &  \qquad  +  \int_0^ t k_2 [\vert \sigma(Q^\varepsilon(s))-\sigma(Q^0(s)) \vert  \vert I^0(s) \vert
+\vert I^\varepsilon(s)-I^0(s) \vert \vert  \sigma(Q^\varepsilon(s) \vert]ds \\
& & \qquad  +  \int_0^ t k_1 [\vert \sigma(Q^\varepsilon(s))-\sigma(Q^0(s)) \vert  \vert S^0(s) \vert
+\vert S^\varepsilon(s)-S^0(s) \vert \vert  \sigma(Q^\varepsilon(s) \vert ds \\ & &\qquad + \varepsilon  \vert J_2(t) \vert ] ds \\ && \quad  \le\int_0^ t
( m +   c_4 C (k_1+k_1 b  e^{-\mu \tau }+ k_2 )  ) \vert  Q^\varepsilon(s) - Q^0(s) \vert ds  \\
& &  \qquad  +  \int_0^ t ( M k_1(1+ b  e^{-\mu \tau })) \vert S^\varepsilon(s)-S^0(s) \vert  ds
 \\ && \qquad
+  \int_0^ t  M k_2 \vert I^\varepsilon(s)-I^0(s) \vert  ds + \varepsilon  \vert J_2(t) \vert
\end{eqnarray*}
and doing the same computations
\begin{eqnarray*}
& & \vert I^\varepsilon(t) - I^0(t)  \vert \le\int_0^ t    \mu \vert  I^\varepsilon(s) - I^0(s) \vert ds   \\ & &\qquad+
\int_0^ t    M k_1(1+e^{-\mu \tau})  \vert S^\varepsilon(s) - S^0(s) \vert ds  \\ & &\qquad + \int_0^ t   k_1 c_4 C (1+e^{-\mu \tau})  \vert Q^\varepsilon(s) - Q^0(s) \vert ds.
\end{eqnarray*}
Thus
$$
\vert Z^\varepsilon(t) - Z^0(t) \vert^ 2 \le c_5 \varepsilon^ 2  ( \vert J_1(t) \vert^ 2 + \vert J_2(t) \vert^ 2 ) + c_6
\int_0^ t \vert Z^\varepsilon(s) - Z^0(s) \vert^ 2 ds. $$
The proof finishes using a Gronwall's lemma type and exponential inequalities for martingales (see the proof of Proposition 3.2 in \cite{BBRT} for the detailed methods).

\hfill$\square$

\end{document}